\documentclass{amsart}

\newtheorem{theorem}{Theorem}[section]

\newtheorem {proposition} [theorem] {Proposition}

\theoremstyle{definition}

\theoremstyle{remark}
\newtheorem{remark}[theorem]{Remark}

\numberwithin{equation}{section}

\begin{document}

\title[Infinitely many solutions for the p-Laplacian in annulus]{Infinitely many solutions for the Dirichlet problem involving the p-Laplacian in annulus}

\author{Anderson L. A. de Araujo}
\address{Departamento de Matem\'atica, Universidade Federal de Vi\c cosa, 36570-000, Vi\c cosa (MG), Brazil}
\email{anderson.araujo@ufv.br}






\begin{abstract}
We present a result of existence of infinitely many solutions for the Dirichlet problem involving the p-Laplacian in annular domains, when $p\leq N$, contouring the failure of compactness of $W^{1,p}(\Omega)$ in $C^0(\overline{\Omega})$ applying a variable change.
\end{abstract}

\maketitle


\par
.

\section{Introduction}
This paper is to investigate the following autonomous Dirichlet problem in annulus
\begin{equation}
\label{P}
\left\{
\begin{array}{cll}
\displaystyle -\Delta_p\,u = f(u) & \textup{in} &\Omega_{a,b},\\
\displaystyle u=0 &\textup{on}& \partial\Omega_{a,b}
\end{array}
\right.
\end{equation} 
 where $\Omega_{a,b}= \{x \in \mathbb{R}^N: a<|x|<b\}$ with $0<a<b$ constants in $\mathbb{R}$, 
\[1<p\leq N,\] 
$\Delta_pu=div(|\nabla\,u|^{p-2}\nabla\,u)$ and $f:\mathbb{R} \to \mathbb{R}$ is a continuous function.

In order to study the solution of \eqref{P}, one can make a standard change of variables. In the $N>p$, be $t=-\frac{A}{r^{(N-p)/(p-1))}}+B$ and $v(t)=u(r)$, where 
\[A=\frac{(ab)^{\frac{N-p}{p-1}}}{b^{\frac{N-p}{p-1}} - a^{\frac{N-p}{p-1}}}\,\,\,\mbox{and}\,\,\, B=\frac{b^{\frac{N-p}{p-1}}}{b^{\frac{N-p}{p-1}} - a^{\frac{N-p}{p-1}}},\] 
 then the problem \eqref{P} transforms into the boundary value problem for the nonautonomous ODE
\begin{equation}
\label{P2}
\left\{
\begin{array}{lll}
\displaystyle (|v'(t)|^{p-2}v'(t))'+q(t)f(v(t)) = 0 & \textup{in} &(0,1),\\
\displaystyle v(0)=v(1)=0, && 
\end{array}
\right.
\end{equation}
where 
\begin{equation}\label{q}
	q(t)=\left(\frac{p-1}{N-p}\right)^{p}\frac{A^{\frac{(p-1)p}{N-p}}}{(B-t)^{\frac{p(N-1)}{N-p}}}.
\end{equation}

In the case $p=N$, one sets $r=a(\frac{b}{a})^t$ and $v(t)=u(r)$, obtaining again the problem \eqref{P2}, now with
\[q(t)=\left[a(\frac{b}{a})^t\ln\frac{b}{a}\right]^{p}.\]
Note that in both cases, the function $q(t)$ is well defined, continuous and bounded between positive constants in the interval $[0,1]$, that is, there exist $q_1,q_0>0$ such that 
\[0<q_0\leq q(t)\leq q_1.\] 

For our purpose, we shall restrict our attention to the ordinary boundary value problem \eqref{P2},
where the function $q(t)$ is continuous and positive on the interval $[0,1]$, while for $f$ we consider the
assumptions below. A weak solution of \eqref{P2} is any $v \in W^{1,p}_0(0,1)$ such that
\[\int_0^1|v'(t)|^{p-2}v'(t)w'(t)dt - \int_0^1q(t)f(v(t))w(t)dt=0,\]
for each $w \in W^{1,p}_0(0,1)$.

We are interested in the existence of infinitely many non-negative weak solutions for
problem \eqref{P}, or equivalently, to the problem \eqref{P2}. Precisely, if $F(\xi)=\int_0^{\xi}f(t)dt$, our aim is to prove the following results:
\begin{theorem}\label{T1}
Assume that $f(0)=0$, $q(t)\geq q_0>0$ in $[0,1]$ and $\inf_{\xi\geq 0}F(\xi)\geq 0$. Moreover, suppose that there exist two sequences $\{a_k\}_{k \in \mathbb{N}}$ and $\{b_k\}_{k \in \mathbb{N}}$ in $]0,+\infty[$, with $a_k<b_k$, $\lim_{k \to +\infty}b_k=+\infty$, such that
\begin{description}
	\item[(i)]$\lim_{k\to +\infty}\frac{b_k}{a_k}=+\infty$;
	\item[(ii)] $\max_{[a_k,b_k]}f\leq 0$ for all $k \in \mathbb{N}$;
	\item[(iii)] Suppose that
	\[\frac{\sigma(p,q_0)}{p(\sup_{t \in (0,1)}dist(t,\{0,1\}))^p}<\limsup_{\xi \to +\infty}\frac{F(\xi)}{\xi^p}<+\infty,\]
	where
	\[\sigma(p,q_0)=\inf_{\mu \in ]0,1[}\frac{1}{q_0\mu(1-\mu)^{p-1}}.\]
\end{description}
Then, problem \eqref{P2} (respectively the problem \eqref{P}) admits an unbounded sequence of non-negative weak solutions in $W^{1,p}_0(0,1)$ (respectively in $W^{1,p}_0(\Omega_{a,b})$).
\end{theorem}

\begin{remark}\label{1.1}
The constant $\sigma(p,q_0)$ is well defined. To see this, let $\sigma:(0,1) \to \mathbb{R}$ be the function defined by $\sigma(x)=\frac{1}{x(1-x)^{p-1}}$. As   $\sigma'(x)=\frac{-(1-x)^{p-2}(1-px)}{x^2(1-x)^{2(p-1)}}$ we have $\sigma'(x)< 0$ if $0<x< \frac{1}{p}$, $\sigma'(\frac{1}{p})=0$ and $\sigma'(x)> 0$ if $\frac{1}{p}< x <1$. Therefore,
\[\sigma(p,q_0)=\inf_{\mu \in ]0,1[}\frac{1}{q_0\mu(1-\mu)^{p-1}} = \frac{1}{q_0p}.\]
\end{remark}

\begin{theorem}\label{T2}
Assume that $f(0)=0$, $q(t)\geq q_0>0$ in $[0,1]$ and $\inf_{\xi\geq 0}F(\xi)\geq 0$. Moreover, suppose that there exist two sequences $\{a_k\}_{k \in \mathbb{N}}$ and $\{b_k\}_{k \in \mathbb{N}}$ in $]0,+\infty[$, with $a_k<b_k$, $\lim_{k \to +\infty}b_k=+\infty$, such that
\begin{description}
	\item[(i)]$\lim_{k\to +\infty}\frac{b_k}{a_k}=+\infty$;
	\item[(ii)] $\max_{[a_k,b_k]}f\leq 0$ for all $k \in \mathbb{N}$;
	\item[(iii)] Suppose that
	\[\frac{\sigma(p,q_0)}{p(\sup_{t \in (0,1)}dist(t,\{0,1\}))^p}<\limsup_{\xi \to 0^+}\frac{F(\xi)}{\xi^p}<+\infty,\]
	where
	\[\sigma(p,q_0)=\inf_{\mu \in ]0,1[}\frac{1}{q_0\mu(1-\mu)^{p-1}}.\]
\end{description}
Then, problem \eqref{P2} (respectively the problem \eqref{P}) admits an sequence of non-negative weak solutions, which strongly converges to $0$ in $W^{1,p}_0(0,1)$ (respectively in $W^{1,p}_0(\Omega_{a,b})$).
\end{theorem}

The statement and the proves of Theorems \ref{T1} and \ref{T2} are very similar the Theorems $1.1$ and $1.2$ of F. Cammaroto, A. Chinnì, B. Di Bella \cite{CCD} which is the principal motivation and reference of this paper. But in \cite{CCD} the authors supposed that $p>N$ and a standard argument, chiefly based on the compact embedding $W^{1,p}(\Omega)$ in $C^0(\overline{\Omega})$ while in this paper we consider $p\leq N$ and soon, the above mentioned immersion is not satisfied. For the sake of completeness, we decided to remake the proves of Theorems \ref{T1} and \ref{T2} to explicit the most important changes in specific arguments compared with \cite{CCD}.

The existence of infinitely many solutions of problem \eqref{P} in general bounded domains $\Omega$ has been studied extensively. Among them, the ones which are closest to the present article are certainly \cite{CCD, OZ, OZ2} and references therein. The approach used in \cite{CCD} is based on a recent variational principle obtained by Ricceri in \cite{Ri}, while in \cite{OZ, OZ2} is based on the method of lower and upper solutions. 

A more general problem that \eqref{P} was studied in Bonanno and Bisci \cite{BB} in general bounded domains $\Omega$, and the authors also supposed that $p>N$ following the same arguments of \cite{CCD}. Also examples and applications are given in comparison with \cite{CCD}. Following ours results, we can prove analogous results to the results of \cite{BB}, in $\Omega_{a,b}$, when $p\leq N$.   

The ours approach is based on the changed of variables, problem \eqref{P2}, and the recent variational principle obtained by Ricceri in \cite{Ri}. The following result is a direct consequence of Theorem $2.5$ of \cite{Ri}.  

\begin{proposition}\label{Prop1}
Let $X$ be a reflexive real Banach space, and let $\Phi$, $\Psi: X \to \mathbb{R}$ be two sequentially
weakly lower semicontinuous and G\^ateaux differentiable functional. Assume
also that $\Psi$ is (strongly) continuous and satisfies $\lim_{\|x\|\to +\infty}\Psi(x) = +\infty$. For each $r>\inf_{X}\Psi$, put 
\[\varphi(r)=\inf_{x \in \Psi^{-1}(]-\infty,r[)}\frac{\Phi(x)-\inf_{\overline{(\Psi^{-1}(]-\infty,r[))_w}}\Phi}{r-\Psi(x)},\] 
where $\overline{(\Psi^{-1}(]-\infty,r[))_w}$ is the closure of $\Psi^{-1}(]-\infty,r[)$ in the weak topology. Fixed $\lambda \in \mathbb{R}$, then
\begin{description}
	\item[(1)] if $\{r_n\}_{ \in \mathbb{N}}$ is a real sequences with $\lim_{n \to +\infty}r_n= +\infty$ such that $\varphi(r_n)<\lambda$, for each $n \in \mathbb{N}$, the following alternative holds: either $\Phi + \lambda\Psi$ has a global minimum, or there
exists a sequence $\{x_n\}$ of critical points of $\Phi + \lambda\Psi$ such that $\lim_{n \to +\infty}\Psi(x_n)= +\infty$.
\item[(2)] if $\{s_n\}_{ \in \mathbb{N}}$ is a real sequences with $\lim_{n \to +\infty}s_n= (\inf_X\Psi)^+$ such that $\varphi(s_n)<\lambda$, for each $n \in \mathbb{N}$, the following alternative holds: either there exists a global minimum of $\Psi$  which is a local minimum of $\Phi + \lambda\Psi$, or there exists a sequence $\{x_n\}$ of pairwise distinct critical points of $\Phi + \lambda\Psi$, with $\lim_{n \to +\infty}\Psi(x_n)=\inf_{X}\Psi$, which weakly converges to a global minimum of $\Psi$.
\end{description}
\end{proposition}

Before stating our main result, we wish to point out that in the sequel $f(x) = 0$ for each $x \in ]-\infty,0[$. We shall consider the Sobolev space $W^{1,p}_0(0,1)$ endowed with the norm
\[\|v\|:=\left(\int_0^1|v'(t)|^pdt\right)^{1/p}.\]
We recall that there exists a constant $c>0$ such that
\begin{equation}\label{immersion}
	\sup_{x \in [0,1]}|v(t)|\leq c\|v\|
\end{equation}
for each $v \in W^{1,p}_0(0,1)$.

\section{Proofs }

\begin{proof}[\textbf{Proof of Theorem \ref{T1}}] Let us apply part ($a$) of Proposition \ref{Prop1}. To this end choose $X=W^{1,p}_0(0,1)$ and for each $v \in X$, put
\[\Phi(v)=-\int_0^1\left(\int_0^{u(t)}q(t)f(s)ds\right)dt\]
and
\[\Psi(v)=\|v\|^p.\]
It is well known that the critical points in $X$ of the functional $\Phi + (1/p)\Psi$ are precisely the weak solutions of problem \eqref{P2}. By the compact embedding of $W^{1,p}_0(0,1)$ in $C([0,1])$, it is not difficult ensures that the functionals $\Phi$ and $\Psi$ are G\^ateaux differentiable and sequentially weakly lower semicontinuous, moreover $\Psi$ is obviously (strong) continuous and coercive.

In our case the function $\varphi$ of Proposition \ref{Prop1} is defined by setting
\[\varphi(r)=\inf_{\|v\|^p<r} \frac{\sup_{\|v\|^p<r}\int_0^1q(t)F(v(t))dt - \int_0^1q(t)F(v(t))dt}{r-\|v\|^p}\]
for each $r \in ]0, +\infty[$.

Now, put $r_k=\left(\frac{b_k}{c}\right)^p$, we wish to prove that $\varphi(r_k)<\frac{1}{p}$ for each $k \in \mathbb{N}$. To this aim, it
suffices to prove that, for each $k \in \mathbb{N}$, there exists a function $v_k \in X$, with $\|v_k\|^p<r_k$, such that
\[\sup_{\|v\|^p<r_k}\int_0^1q(t)F(v(t))dt - \int_0^1q(t)F(v_k(t))dt<\frac{1}{p}(r_k-\|v_k\|^p).\]
From ($iii$) we can choose a constant $h$ such that
\begin{equation}\label{est1}
	\frac{\sigma(p,q_0)}{p(\sup_{t \in (0,1)}dist(t,\{0,1\}))^p}<h<\limsup_{\xi \to +\infty}\frac{F(\xi)}{\xi^p}
\end{equation}
and so there exists $t_0 \in (0,1)$ such that $\left(\frac{\sigma(p,q_0)}{ph}\right)^{1/p}< dist(t_0, \{0,1\})$. Therefore, we can fix $\gamma$
satisfying
\[\left(\frac{\sigma(p,q_0)}{ph}\right)^{1/p}<\gamma< dist(t_0, \{0,1\}).\]
Observe that by ($ii$), 
\begin{equation}\label{max}
	\max_{[0,a_k]}F=\max_{[0,b_k]}F.
\end{equation}
 Now, fix $k \in \mathbb{N}$ and consider the function $v_k \in X$ defined by setting
\[
v_k(t)=\left\{
\begin{array}{lll}
\displaystyle 0 & \textup{if} &t \in (0,1)\backslash(t_0-\gamma,t_0+\gamma),\\
\displaystyle \xi_k & \textup{if} &t \in (t_0-\frac{\gamma}{2},t_0+\frac{\gamma}{2}),\\
\displaystyle \frac{2\xi_k}{\gamma}(\gamma - |t-t_0|) & \textup{if} &t \in (t_0-\gamma,t_0+\gamma)\backslash(t_0-\frac{\gamma}{2},t_0+\frac{\gamma}{2})\\
\end{array}
\right.
\]
with $\xi_k \in ]0,a_k]$ such that 
\begin{equation}\label{max2}
	F(\xi_k)=\max_{\eta \in [0,a_k]}F(\eta).
\end{equation}

In view of ($i$), we can choose $k \in \mathbb{N}$ so that
\begin{equation}\label{est2}
\frac{b_k}{a_k}>\frac{2c}{\gamma}\left(\gamma\right)^{1/p}	
\end{equation}
fro all $k > \overline{k}$.

Fix $k > \overline{k}$. By \eqref{est2} we have
\[
\begin{array}{rcl}
\|v_k\|^p&=&\displaystyle \int_0^1|v_k'(t)|^pdt=\int_{(t_0-\gamma,t_0+\gamma)\backslash(t_0-\frac{\gamma}{2},t_0+\frac{\gamma}{2})}\frac{2^p\xi_k^p}{\gamma^p}dt\\
&=& \displaystyle\frac{2^p\xi_k^p}{\gamma^p}(|(t_0-\gamma,t_0+\gamma)|-|(t_0-\frac{\gamma}{2},t_0+\frac{\gamma}{2})|)=\frac{2^p\xi_k^p}{\gamma^p}(2\gamma - \gamma)\\
&=& \displaystyle\frac{2^p\xi_k^p}{\gamma^{p-1}} <r_k.
\end{array}
\]
Because of \eqref{immersion}, \eqref{max} and \eqref{max2}, for each $v \in X$, fulfilling $\|v\|^p\leq r_k$, one has
\[\sup_{x \in [0,1]}|v(t)|\leq cr_k^{1/p}\]
as well as
\[F(v(t))\leq \max_{\eta \in [0,cr_k^p]}F(\eta)=F(\xi_k)\,\,\,\mbox{in}\,\,\,(0,1).\]
Next, since $\lim_{k \to +\infty}\frac{r_k}{\xi_k^p}=+\infty$, there exists $k_0 \in \mathbb{N}$ such that
\begin{equation}\label{est3}
\frac{r_k}{\xi_k^p}>p\limsup_{\xi \to +\infty}\frac{F(\xi)}{\xi^p}\left(\int_0^1q(t)dt - \int_{t_0-\frac{\gamma}{2}}^{t_0+\frac{\gamma}{2}}q(t)dt\right) + \frac{2^p}{\gamma^{p-1}}
\end{equation}
for all $k>k_0$. Hence, using \eqref{est3}, we get
\[
\begin{array}{r}
\displaystyle \sup_{\|v\|^p<r_k}\int_0^1q(t)F(v(t))dt - \int_0^1q(t)F(v_k(t))dt \\
 \displaystyle \leq F(\xi_k)\int_0^1q(t)dt - \int_{t_0-\frac{\gamma}{2}}^{t_0+\frac{\gamma}{2}}q(t)F(\xi_k)dt\\
 \displaystyle \leq F(\xi_k)\left(\int_0^1q(t)dt - \int_{t_0-\frac{\gamma}{2}}^{t_0+\frac{\gamma}{2}}q(t)dt\right)\\
 \displaystyle <\frac{1}{p}(r_k-\|v_k\|^p)
\end{array}
\]
for each $k>k^*\geq k_0$.

Since $\lim_{k \to +\infty}r_k=+\infty$, the previous inequality assures that the conclusion ($1$) of Proposition \ref{Prop1} can be used and either the functional $\Phi + (1/p)\Psi$ has a global minimum, or there exists a sequence $\{v_k\}_{k\in \mathbb{N}}$ of solutions of problem \eqref{P2} such that $\lim_{k \to +\infty}\|v_k\|=+\infty$.

The other step is to verify that the functional $\Phi + (1/p)\Psi$ has no global minimum. Taking into account \eqref{est1}, one has, for each $k\in \mathbb{N}$.
\[\sup_{\eta\geq k}\frac{F(\eta)}{\eta^p}>h\]
and so there exists $\eta_k\geq k$ such that
\[\frac{F(\eta_k)}{\eta_k^p}>h.\]
Now, if we consider a function $w_k \in X$ defined by setting
\[
w_k(t)=\left\{
\begin{array}{lll}
\displaystyle 0 & \textup{if} &x \in (0,1)\backslash(t_0-\gamma,t_0+\gamma),\\
\displaystyle \eta_k & \textup{if} &x \in (t_0-\overline{\mu}\gamma,t_0+\overline{\mu}\gamma),\\
\displaystyle \frac{\eta_k}{\gamma(1-\overline{\mu})}(\gamma - |t-t_0|) & \textup{if} &x \in (t_0-\gamma,t_0+\gamma)\backslash(t_0-\overline{\mu}\gamma,t_0+\overline{\mu}\gamma)\\
\end{array}
\right.
\]
where $\overline{\mu} \in ]0, 1[$ is such that
\[\sigma(p,q_0)=\frac{1}{q_0\overline{\mu}(1-\overline{\mu})^{p-1}},\]
note that $\overline{\mu}$ is well defined by Remark \ref{1.1}. We have
\[
\begin{array}{l}
\Phi(w_k) + \frac{1}{p}\Psi(w_k)\\
=\displaystyle -\int_0^1q(t)F(w_k(t))dt + \frac{1}{p}\|w_k\|^p\\
\leq \displaystyle -\int_{t_0-\overline{\mu}\gamma}^{t_0+\overline{\mu}\gamma}q(t)F(\eta_k)dt + \frac{1}{p}\frac{\eta_k^p}{\gamma^p(1-\overline{\mu})^p}(2\gamma - 2\overline{\mu}\gamma)\\
\leq \displaystyle -2\overline{\mu}\gamma\,q_0F(\eta_k) + \frac{2\gamma}{p}\frac{\eta_k^p}{\gamma^p(1-\overline{\mu})^{p-1}}\\
=\displaystyle 2\overline{\mu}\gamma\,q_0(\frac{1}{p}\frac{\eta_k^p}{\gamma^pq_0\overline{\mu}(1-\overline{\mu})^{p-1}} - F(\eta_k))\\
<\displaystyle 2\overline{\mu}\gamma\,q_0(\frac{1}{p}\frac{\eta_k^p}{\gamma^pq_0\overline{\mu}(1-\overline{\mu})^{p-1}} - h\eta_k^p)=2\overline{\mu}\gamma\,q_0\eta_k^p\left(\frac{\sigma(p,q_0)}{p\gamma^p} - h\right).\\
\end{array}
\]
Since $h>\frac{\sigma(p,q_0)}{p\gamma^p}$, we conclude that $\lim_{k \to +\infty}\eta_k^p\left(\frac{\sigma(p,q_0)}{p\gamma^p} - h\right)=-\infty$ and so the previous inequality
shows that the functional $\Phi + (1/p)\Psi$ is not bounded from below and then it has no global minimum.

Therefore, Proposition \ref{Prop1} assures that there is a sequence $\{v_k\}_{k \in \mathbb{N}} \subset X$ of critical points of $\Phi + (1/p)\Psi$ such that $\lim_{n \to +\infty}\|v_n\|= +\infty$. As previously observed, every function $v_k$ is a weak solution of \eqref{P2}.

Finally, we claim that each weak solution of problem \eqref{P2} is non-negative in $(0,1)$. Assume the contrary. Let $v$ be a weak solution of \eqref{P2} such that the set $A=\{x \in (0,1): v(x)<0\}$ is non-empty. By the continuity of $v$, $A$ is open and so $v|_{A} \in W^{1,p}_0(A)$. Then, for each $w \in W^{1,p}_0(A)$,
\[\int_A|v'(t)|^{p-2}v'(t)w'(t)dt - \int_Aq(t)f(v(t))w(t)dt=0.\]
The assumptions on $f$ imply
\[\int_A|v'(t)|^{p-2}v'(t)w'(t)dt=0\]
for each $w \in W^{1,p}_0(A)$ and so, in particular, one has 
\[\int_A|v'(t)|^{p}dt=0,\]
an absurd. This completes the proof.
\end{proof}

\begin{proof}[\textbf{Proof of Theorem \ref{T2}}] 

We take $X, \Phi, \Psi$ as in the proof of Theorem \ref{T1}. In a similar way we prove that $\varphi(s_k)< \frac{1}{p}$ for each $k \in \mathbb{N}$, with $s_k =\left(\frac{b_k}{c}\right)^p$. Now, fixed $h$ such that
\[	\frac{\sigma(p,q_0)}{p(\sup_{t \in (0,1)}dist(t,\{0,1\}))^p}<h<\limsup_{\xi \to 0^+}\frac{F(\xi)}{\xi^p}.\]
For each $k \in \mathbb{N}$, there exists $\eta_k\leq \frac{1}{k}$ such that $\frac{F(\eta_k)}{\eta_k^p}>h$. If we take $w_k$ as in the proof
of Theorem \ref{T1}, of course the sequence $\{w_k\}$ strongly converges to 0 in $X$ and $\Phi(w_k) + (1/p)\Psi(w_k)<0$ for all $k \in \mathbb{N}$. Since $\Phi(0) + (1/p)\Psi(0)=0$, this means that $0$ is not a local minimum of $\Phi + (1/p)\Psi$. Then, since $0$ is the only a global minimum of $\Psi$, the part ($2$) of Proposition \ref{Prop1} ensures that there exists a sequence $\{v_k\}_{k \in \mathbb{N}} \subset X$ of critical points of $\Phi + (1/p)\Psi$ such that $\lim_{k \to +\infty}\|v_k\|=0$ and this completes the proof.
\end{proof}

\section{Conclusions}

We can also to investigate the following nonautonomous Dirichlet problem in an annular domain
\begin{equation}
\label{P3}
\left\{
\begin{array}{cll}
\displaystyle -\Delta_p\,u = g(|x|)f(u) & \textup{in} &\Omega_{a,b},\\
\displaystyle u=0 &\textup{on}& \partial\Omega_{a,b}
\end{array}
\right.
\end{equation} 
 where $\Omega_{a,b}= \{x \in \mathbb{R}^N: a<|x|<b\}$ with $0<a<b$ constants in $\mathbb{R}$, $1<p\leq N$, 
$\Delta_pu=div(|\nabla\,u|^{p-2}\nabla\,u)$, $f:\mathbb{R} \to \mathbb{R}$ is a continuous function and $g:[a,b] \to \mathbb{R}$  is a continuous function such that $g(s)\geq g_0>0$.

In order to study the solution of \eqref{P3}, one can make a standard change of variables, see Liu and Yang \cite{LY}. In the $N\geq p+1$, if $t=-\frac{A}{r^{(N-p)/(p-1))}}+B$ and $v(t)=u(r)$, where 
\[A=\frac{(ab)^{\frac{N-p}{p-1}}}{b^{\frac{N-p}{p-1}} - a^{\frac{N-p}{p-1}}}\,\,\,\mbox{and}\,\,\, B=\frac{b^{\frac{N-p}{p-1}}}{b^{\frac{N-p}{p-1}} - a^{\frac{N-p}{p-1}}},\] 
 then the problem \eqref{P} transforms into the boundary value problem for the nonautonomous ODE
\begin{equation}
\label{P4}
\left\{
\begin{array}{lll}
\displaystyle (|v'(t)|^{p-2}v'(t))'+h(t)k(t)f(v(t)) = 0 & \textup{in} &(0,1),\\
\displaystyle v(0)=v(1)=0, && 
\end{array}
\right.
\end{equation}
where 
\[h(t)=g\left(\left(\frac{A}{B-t}\right)^{\frac{p-1}{N-p}}\right)\]
and
\[k(t)=\left(\frac{p-1}{N-p}\right)^{p}\frac{A^{\frac{(p-1)p}{N-p}}}{(B-t)^{\frac{p(N-1)}{N-p}}}.\]
In the case $p=N$, one sets $r=b\left(\frac{a}{b}\right)^t$ and $v(t)=u(r)$, obtaining again the problem \eqref{P2}, now with
\[h(t)=g\left(b\left(\frac{a}{b}\right)^t\right) \,\,\mbox{and} \,\, k(t)=\left[b\left(\frac{a}{b}\right)^t\left(\ln\frac{b}{a}\right)^{-1}\right]^{p-1}.\]
Note that in both cases, the function $h(t)$ and $q(t)$ are well defined, continuous and bounded between positive constants in the interval $[0,1]$, that is, there exists $q_0,q_1,q_2, q_3>0$ such that 
\[0<q_0\leq h(t)\leq q_1 \,\,\mbox{and} \,\, 0<q_2\leq k(t)\leq q_3.\] 
Then, in \eqref{P3} put $q(t)=h(t)k(t)$ and we can get analogous results as in Theorems \ref{T1} and \ref{T2}.

\bibliographystyle{amsplain}

\end{document}